\documentclass[11pt]{amsart}
\usepackage{amscd,verbatim}
\numberwithin{equation}{section}
\newtheorem{Theorem}{Theorem}[section]
\newtheorem{Lemma}{Lemma}[section]
\newtheorem{Corollary}{Corollary}[section]

\theoremstyle{definition}
\newtheorem{Definition}{Definition}[section]

\theoremstyle{remark}
\newtheorem{Remark}{Remark}

\newcommand{\bfR}{\mathbb{ R}}
\newcommand{\bfC}{\mathbb{ C}}

\newcommand{\bfH}{\mathbb{ H}}

\renewcommand{\Im}{{\rm Im}}
\renewcommand{\d}{{\rm d}}

\title[Bonnet pairs]{Bonnet pairs and isothermic surfaces}
\author[G. Kamberov]{George Kamberov}
\address{Department of Mathematics \\
University of Massachusetts \\
Amherst, MA 01003}
\email[George Kamberov]{george@gang.umass.edu}

\author[F. Pedit]{Franz Pedit}
\address{Department of Mathematics \\
University of Massachusetts \\
Amherst, MA 01003}
\email[Franz Pedit]{franz@gang.umass.edu}

\author[U. Pinkall]{Ulrich Pinkall}
\address{Fachbereich Mathematik, MA 8-3 \\
Technische Universit\"at Berlin \\
Strasse des 17. Juni 136\\ 10622 Berlin, Germany}
\email[Ulrich Pinkall]{pinkall@sfb288.math.tu-berlin.de}
\thanks{First and second author were partially supported by NSF grant DMS-9205293.
Third author was supported by a Five Colleges Distinguished Visiting
Professorship and NSF grant DMS-9205293.}
\subjclass{}
\keywords{}

\date{}

\begin{document}

\maketitle
\section{Introduction}
A classical question in surface theory is which data are sufficient
to describe a surface in space up to rigid motions. Bonnet 
suggested that mean curvature and metric should suffice to determine
the surface generically. The local theory was developed by Bonnet
\cite{Bon}, Cartan \cite{Car} and Chern \cite{Che} who showed the existence
of various $1$-parameter families of Bonnet surfaces, i.e., surfaces with the
same induced metric and mean curvature. A comprehensive study of this problem
and its relationship to the Painleve equations has recently been completed
by Bobenko and Eitner \cite{BobEig}. On the other hand, Lawson and Tribuzy
\cite{LawTri} have shown that for embedded compact surfaces
there are at most two surfaces to a given metric and mean curvature.
Moreover, uniqueness can be established under various global
assumptions \cite{Kam}.
Up to date it is unknown whether such compact Bonnet pairs exist and if, how to 
construct them.

In this note we classify {\em all} Bonnet pairs on a simply connected domain.
Our main intent was to apply what we call a 
{\em quaternionic function theory} to a concrete problem in differential geometry.
In the first section we develop the necessary formulas and explain how
usual Riemann surface theory can be viewed as a special case of our
extended function theory. The ideas are simple: conformal immersions into
quaternions or imaginary quaternions take the place of chart maps for a Riemann
surface. Starting from a reference immersion we construct all conformal
immersions of a given Riemann surface (up to translational periods) by
{\em spin transformations} (Definition~\ref{def:spin}).With this viewpoint in mind
we discuss in the second section how to construct all
Bonnet pairs on a simply connected domain from isothermic surfaces
and vice versa.  Isothermic surfaces are solutions to a certain
soliton equation \cite{BurJerPedPin} thus a simple dimension count
tells us that we obtain Bonnet pairs which are not
part of any of the families. The corresponcence between Bonnet pairs
and isothermic surfaces is explicit and to each isothermic surface
we obtain a 4-parameter family of Bonnet pairs.
 
\section{Surface theory revisited}
Let $M$ be a Riemann surface and $f:M\to\bfR^3$ an immersion. We always regard
$\bfR^4=\bfH$ as quaternions and $\bfR^3=\Im\bfH$ as purely imaginary quaternions.
Note that for $a,b\in\Im\bfH$ we have
\[
ab=-<a,b>+a\times b\,.
\]
There is a natural wedge product over $\bfH$-valued $1$-forms given by
\begin{equation}\label{eqn:wedge}
\alpha\wedge\beta(X,Y)=\alpha(X)\beta(Y)-\alpha(Y)\beta(X)
\end{equation}
which satisfies the obvious identities
\begin{gather*}
\overline{\alpha\wedge\beta}=-\bar{\beta}\wedge\bar{\alpha}\\
\alpha\wedge h\beta=\alpha h\wedge\beta\\
\d(h\alpha)=d h\wedge\d\alpha+h\d\alpha\\
\d(\alpha h)=\d\alpha h-\alpha\wedge dh
\end{gather*}
for an $\bfH$-valued function $h:M\to\bfH$. Usually we will identify $2$-forms on $M$
with their quadratic forms via
\[
\omega(X)=\omega(X,JX)
\]
where $\omega\in\Omega^2(M,\bfH)$ and $J:TM\to TM$ is the complex structure of $M$. Thus the
only objects appearing in our theory will be $\bfH$-valued functions on $M$ and $TM$. A
key formula is
\begin{equation}\label{eqn:key}
\alpha\wedge\beta=\alpha(*\beta)-(*\alpha)\beta
\end{equation}
where $*\alpha=\alpha\circ J$ is minus the usual Hodge star operator. With this in place
let us develop the necessary equations of surface theory as needed for our purposes.
\begin{Lemma}
$f:M\to\bfR^3$ is a conformal immersion if and only if  there exists $N:M\to\bfH$ such that
\begin{equation}\label{eqn:conformal}
*d f=Nd f\,.
\end{equation}
If \eqref{eqn:conformal} holds then $N:M\to S^2\subset\Im\bfH$ is the (oriented) unit normal
field (Gauss map) to $f$.
\end{Lemma}
\begin{proof}
Since $*^2=-1$ and $d f$ is pointwise injective we have $N^2=-1$. This implies $|N|^2=1$
and thus $\bar{N}=-N$. Conjugating \eqref{eqn:conformal} gives 
$Nd f=-d f N$ which says that $<N,d f>=0$. Thus $N:M\to S^2\subset\Im\bfH$ is a unit
normal field to $f$ which defines the orientation since $*d f=d f\circ J$.
\end{proof}
\begin{Remark}
Equation \eqref{eqn:conformal} is a generalization of the Cauchy-Riemann equation:
to $h:M\to\bfC$ let $f=hj:M\to\Im\bfH$ with $N=i$. Then \eqref{eqn:conformal} becomes
\[
*d h= id h
\]
which is the Cauchy-Riemann equation for $h$. In this sense surface theory in $\Im\bfH$
can be viewed as a generalization of complex function theory to a quaternionic
function theory. Holomorphic chart maps are special cases of conformal immersions.
\end{Remark}
Any $\bfR^3$-valued $1$-form $\alpha$ on $M$ can be decomposed into its conformal
and anti-conformal parts:
\begin{subequations}\label{eqn:parts}
\begin{gather}
\alpha=\alpha_{+}+\alpha_{-}\,,\qquad *\alpha_{\pm}=\pm N\alpha_{\pm}\\
\intertext{where}
\alpha_{\pm}=\frac{1}{2}(\alpha\pm N*\alpha)\,.
\end{gather}
\end{subequations}
Thus, using \eqref{eqn:conformal},
\[
d *d f=d N\wedge d f=d N*d f-*d Nd f=(d NN-*d N)d f=2(d N)_{+}Nd f
\]
which is $\Im\bfH$-valued. This means
\[
<(d N)_{+}, Nd f>=0
\]
or, since $d N_{+}$ is also tangential,
\begin{equation}\label{eqn:mean}
d N_{+}=Hd f
\end{equation}
for some real valued function $H:M\to\bfR$. Clearly, $H$ is the mean curvature of $f$
w.r.t. the induced metric $|df|^2$.
Inserting \eqref{eqn:mean} we obtain
\begin{equation}\label{eqn:Weingarten}
\d N=Hd f+\omega\,,\qquad \omega=(d N)_{-}\,,\qquad *\omega=-N\omega
\end{equation}
and thus
\begin{equation}\label{eqn:Laplace}
d *d f=2HN|d f|^2\,.
\end{equation}
Differentiating \eqref{eqn:Weingarten} we get the Codazzi equation
\begin{equation}\label{eqn:Codazzi}
\d\omega=(*d H-d HN)d f\,.
\end{equation}
\begin{Definition}\label{def:spin}
Two conformal immersions $f,\tilde{f}:M\to\bfR^3$ are called {\em spin-equivalent}
if there exists $\lambda:M\to\bfH_{*}$ so that 
\begin{equation}\label{eqn:spin}
d\tilde{f}=\bar{\lambda}d f\lambda\,.
\end{equation}
We call $\tilde{f}$ a {\em spin-transform} of $f$.
\end{Definition}
\begin{Remark}
\begin{enumerate}
\item If $\pi_{1}(M)=0$ then any two conformal immersions are spin-equivalent: since
$d\tilde{f}$ and $d f$ are conformal $1$-forms pointwise they can be mapped into each
other by a rotation and scaling. But this is exactly what \eqref{eqn:spin} means.
Notice that $\lambda$ is unique up to sign.
\item If $\pi_{1}(M)\neq 0$ then $f$ and $\tilde{f}$ are spin-equivalent if and only if
they belong to the same regular homotopy class \cite{Pin}.
\item $\lambda$ is constant if and only if $\tilde{f}$ is obtained from $f$ by
an Euclidean motion and scaling. If $|\lambda|=1$ then $\tilde{f}$ and $f$ are
congruent immersions.
\item Spin-equivalence clearly is an equivalence relation.
\end{enumerate}
\end{Remark}
To obtain new conformal immersions (locally) from a given {\em reference} immersion
$f:M\to\bfR^3$ via spin-transformations we only have to solve
\[
0=d(\bar{\lambda}d f\lambda)=d\bar{\lambda}\wedge d f\lambda-
  \bar{\lambda} d f\wedge d\lambda=-2\Im(\bar{\lambda} d f\wedge d\lambda)\,.
\]
But this says that $\bar{\lambda} d f\wedge d\lambda$ is real valued and hence there
exists a unique real valued function $\rho:M\to\bfR$ with
\[
 d f\wedge d\lambda=-\rho|df|^2\lambda=\rho df^2\lambda\,.
\]
Using \eqref{eqn:key}, \eqref{eqn:conformal} and dividing by $df$ we obtain
the integrability equation
\begin{equation}\label{eqn:Dirac}
*d\lambda+Nd\lambda=\rho df\lambda\,.
\end{equation}
This is worth formulating as a
\begin{Lemma}
If $f,\tilde{f}:M\to\bfR^3$ are spin-equivalent via $d\tilde{f}=\bar{\lambda}df\lambda$
then $\lambda:M\to\bfH_{*}$ satisfies \eqref{eqn:Dirac}.

Conversely, if $\pi_{1}(M)=0$ then to a given conformal immersion $f:M\to\bfR^3$ nowhere vanishing
solutions to \eqref{eqn:Dirac} yield all conformal immersions $\tilde{f}:M\to\bfR^3$ via
spin-transformations $d\tilde{f}=\bar{\lambda}d f\lambda$.
\end{Lemma}
We now relate elementary geometric data of spin-equivalent immersions.
\begin{Lemma}
Let $f,\tilde{f}:M\to\bfR^3$ be spin-equivalent via $d\tilde{f}=\bar{\lambda}d f\lambda$. Then
\begin{enumerate}
\item $\tilde{N}=\lambda^{-1}N\lambda$ is the oriented normal to $\tilde{f}$,
\item $|d\tilde{f}|^2=|\lambda|^4|df|^2$,
\item $\tilde{H}=\frac{H+\rho}{|\lambda|^2}$ where $\rho:M\to\bfR$ is given 
by \eqref{eqn:Dirac}.
\end{enumerate}
\end{Lemma}
\begin{proof}
(i) and (ii) are obvious. To see (iii) we use \eqref{eqn:Weingarten}:
\begin{align*}
\tilde{H}d\tilde{f}+\tilde{\omega}=d\tilde{N}=
\frac{H}{|\lambda|^2}d\tilde{f}+
\lambda^{-1}(\omega+Nd\lambda\lambda^{-1}-d\lambda\lambda^{-1}N)\lambda\,,
\end{align*}
hence, computing the conformal part \eqref{eqn:parts} of the latter
and inserting  \eqref{eqn:Dirac}, we get
\begin{gather*}
(\lambda^{-1}(\omega+Nd\lambda\lambda^{-1}-d\lambda\lambda^{-1}N)\lambda)_{+}=
\lambda^{-1}(\omega+Nd\lambda\lambda^{-1}-d\lambda\lambda^{-1}N)_{+}\lambda=\\
\lambda^{-1}\rho df\lambda=\frac{\rho d\tilde{f}}{|\lambda|^2}
\end{gather*}
and thus
\[
\tilde{H}=\frac{H+\rho}{|\lambda|^2}\,.
\]
\end{proof}
\begin{Corollary}\label{cor:Hdf}
Let $f,\tilde{f}:M\to\bfR^3$ be spin-equivalent via $d\tilde{f}=\bar{\lambda}d f\lambda$.
Then the following are equivalent:
\begin{enumerate}
\item
$\tilde{H}|d\tilde{f}|=H|df|$
\item
$df\wedge\d\lambda=0$ which is the same as $*d\lambda+N d\lambda=0$.
\end{enumerate}
\end{Corollary}
\section{Isothermic surfaces and Bonnet pairs}
We now give the local classification of all Bonnett pairs and explain how they are
obtained from
isothermic surfaces. 
\begin{Definition}\label{def:isothermic}
A conformal immersion $f:M\to\bfR^3$ is {\em isothermic} if there exists an non-zero, 
$\Im\bfH$-valued, closed
and anti-conformal $1$-form $\tau$$\in\Omega^1(M,\Im\bfH)$, i.e., 
\[
\d\tau=0\,,\qquad \tau\neq 0\,,\qquad *\tau+N\tau=0\,\,\text{which is equivalent to}\,\, 
df\wedge\tau=0\,.
\]
Note that this is a Moebius invariant condition. 
\end{Definition}
\begin{Remark}
\begin{enumerate}
\item The classical notion of an isothermic surface is that $f:M\to\bfR^3$ admits
conformal curvature line parameters (away from umbilic points). If $(x,y)$ are such
parameters then one easily sees that the $1$-form
\[
\tau=\frac{1}{|f_x|^2}(f_x dx-f_y dy)=-f_x^{-1}dx+f_y^{-1}dy
\]
is closed and anti-conformal. Thus, classically isothermic surfaces are isothermic
in our notion.
\item Since $\d\tau=0$ locally $\tau=d f^{*}$ and with $N^{*}=-N$ one sees that $f^{*}$
is a conformal immersion (a dual isothermic surface) away from zeros of $\tau$.
\item
Examples of isothermic surfaces include quadrics, surfaces of revolution and 
surfaces of constant mean curvature (and their Moebius transforms). The latter case
follows at once from \eqref{eqn:Codazzi} and \eqref{eqn:Weingarten}: $\d\omega=0$
and $*\omega+N\omega=0$ so that we may put $\tau=\omega$. For minimal surfaces $(H=0)$
we get $dN=\omega$ so that the Gauss map $N$ is a dual surface.
\item
Recently it has been shown \cite{BurJerPedPin, BobPin, JerPed} that isothermic surfaces can be 
obtained from an infinite dimensional integrable system (soliton equation).
This provides an infinte dimesnional space of examples of local isothermic surfaces.
\end{enumerate}
\end{Remark}
\begin{Definition}
Two conformal immersions $f,\tilde{f}:M\to\bfR^3$ form a {\em Bonnet pair} if they
induce the same metric $|d\tilde{f}|^2=|df|^2$ and have the same mean curvature
$\tilde{H}=H$, but are not congruent.
\end{Definition}
Finally we arrive at the local classification of {\em all} Bonnet pairs. We assume
$\pi_{1}(M)=0$.
\begin{Theorem}\label{thm:class}
Let $f:M\to\bfR^3$ be isothermic with dual surface $f^{*}:M\to\bfR^3$.
Choose $\epsilon\in\bfR_{*}$, $a\in\Im\bfH$ and let 
$\lambda_{\pm}=\pm\epsilon+f^{*}+a$. Then the spin transforms $f_{\pm}:M\to\bfR^3$
given by $df_{\pm}=\overline{\lambda_{\pm}}df\lambda_{\pm}$ form a Bonnet pair.

Conversely, every Bonnet pair arises from a $3$-parameter family (determined up to
scalings) of isothermic surfaces
where the three parameters account for Euclidean rotations of the partners in the 
Bonnet pair w.r.t. each other.
\end{Theorem}
\begin{proof}
From Definition~\ref{def:isothermic} we have  $df\wedge df^{*}=0$ and thus also 
$df\wedge d\lambda_{\pm}=0$, i.e., $*d\lambda_{\pm}+Nd\lambda_{\pm}=0$, 
which is \eqref{eqn:Dirac}
with $\rho_{\pm}=0$. Moreover, $|\lambda_{+}|=|\lambda_{-}|$ so that 
$|df_{+}|^2=|df_{-}|^2$ and $H_{+}=H_{-}$. Congruence of $f_{+}$ and $f_{-}$
means that $\lambda_{+}^{-1}\lambda_{-}$ is constant which would imply that $f^{*}$
is constant. Hence $f_{\pm}$ form a Bonnet pair.

To show the converse assume that we have a Bonnet pair $f_{\pm}$. Then
\[
df_{+}=\bar{\lambda}df_{-}\lambda
\]
for some $\lambda:M\to\bfH$. But $f_{+}$ and $f_{-}$ induce the same
metric and same mean curvature so that $|\lambda|=1$ and by 
Corollary~\ref{cor:Hdf}
$df_{-}\wedge\d\lambda=0$. Solving for $f^{*}$ in terms of 
$\lambda_{+}^{-1}\lambda_{-}$ in the first part of the Theorem (with $a=0$) 
suggests to define
\[
f^{*}=(\lambda-1)^{-1}+1/2\,.
\]
To insure that $\lambda-1$ vanishes nowhere we may choose to multiply $\lambda$
on the right by a unit quaternion (which induces a rotation of the partners
of the Bonnet pair w.r.t. each other).
Since $|\lambda|=1$ one easily checks that $f^{*}$ is purely imaginary. Putting
\[
df=\overline{(\lambda-1)}df_{-}(\lambda-1)
\]
we see that $df\wedge df^{*}=0$ so that $f:M\to\bfR^3$ is isothermic with dual 
surface $f^{*}$.
Now one easily checks that $f_{\pm}$ are spin transforms of $f$ via
$\lambda_{\pm}=\pm 1/2 + f^{*}$.
\end{proof}
\begin{Remark}
\begin{enumerate}
\item
Choices of real $\epsilon$ and imaginary quaternionic $a$ give a $4$-parameter family
of Bonnet pairs to each isothermic surface. For $\epsilon=0$ the pair degenerates to a
$3$-paramter family of 
isothermic surfaces $d\tilde{f}=\overline{(f^{*}+a)}df(f^{*}+a)$. Compactifying 
this $\bfR^4$ to $S^4$ adds $f$ (up to congruence) at infinity. 
Thus the above construction provides an
$S^4$-worth of Bonnet pairs with an equatorial $S^3$ of isothermic surfaces.
\item
It is known that there exists at most two compact embedded surfaces in $\bfR^3$
with prescribed mean curvature and metric \cite{LawTri}. Understanding the period problem
in our local classification could provide an affirmative solution to
the question of whether a given compact embedded surface allows a Bonnet partner.
\end{enumerate}
\end{Remark}

\bibliographystyle{amsplain}

\begin{thebibliography}{10}

\bibitem{BurJerPedPin}F. Burstall, U. Hertrich-Jeromin, F. Pedit and U. Pinkall.
	{\em Curved flats and Isothermic surfaces\/}. To appear in Math. Z.
\bibitem{BobEig} A. Bobenko and U. Eitner. {\em Bonnet surfaces and Painleve equations\/}.
SFB preprint  151, 1994.
\bibitem{BobPin} A. Bobenko and U. Pinkall. {\em Discrete isothermic surfaces\/}.
J. reine angew. Math. {\bf 457}, (1996), 187-208.
\bibitem{Car} E. Cartan. {\em Sur les couples de surfaces applicables avec conservation
des courbures principales\/}. Bull. Sc. Math. {\bf 66}, (1942).
\bibitem{Che} S.S. Chern. {\em Deformations of surfaces preserving principal curvatures\/}.
Differential Geometry and Complex Analysis: a volume dedicated to the memory of H.E. Rauch.
Springer Verlag 1985, 155-163.
\bibitem{Bon} O. Bonnet. {\em Memoire sur la theorie des surfaces applicables\/}.
J. Ec. Polyt. {\bf 66}, (1842).
\bibitem{JerPed} U. Hertrich-Jeromin and F. Pedit. {\em Remarks on the Darboux transform of isothermic
surfaces\/}. GANG preprint 1996.
\bibitem{Kam} G. Kamberov. {\em Recovering the shape of a surface from the mean curvature\/}.
Dynamical Systems and Differential equations. M. Dekker 1996. Submitted.
\bibitem{LawTri}B. Lawson and R. Tribuzy. {\em On the mean curvature function for compact
surfaces\/}. J. Diff. Geom. {\bf 3}, (1987).
\bibitem{Pin} U. Pinkall. {\em Regular homotopy classes of immersed surfaces\/}. Topology {\bf 24},
(1985), 412-434.

\end{thebibliography}

\ifx\undefined\bysame
\newcommand{\bysame}{\leavevmode\hbox to3em{\hrulefill}\,}
\fi

\end{document}